\documentclass[12pt]{amsart}

\usepackage{amsmath,xspace,amssymb,mathrsfs}
\usepackage{color}

\input xy
\xyoption{all}
\xyoption{2cell}
\UseAllTwocells
\CompileMatrices

\renewcommand{\phi}{\varphi}

\newcommand{\Ker}{\operatorname{Ker}}

\newcommand{\Supp}{\operatorname{Supp}}

\newtheorem{proposition}{Proposition}[section]
\newtheorem{lemma}[proposition]{Lemma} 
\newtheorem{corollary}[proposition]{Corollary}
\newtheorem{theorem}[proposition]{Theorem}

\newtheorem{prop-def}[proposition]{Proposition and definition}

\theoremstyle{definition}

\newtheorem{remark}[proposition]{Remark}

\begin{document}

\title{A fresh look into monoid rings and formal power series rings}

\author{Abolfazl Tarizadeh}
\address{ Department of Mathematics, Faculty of Basic Sciences, University of Maragheh \\
P. O. Box 55136-553, Maragheh, Iran.
 }
\email{ebulfez1978@gmail.com}


\begin{abstract} In this article, the ring of polynomials is studied in a systematic way through the theory of monoid rings. As a consequence, this study provides natural and canonical approaches in order to find easy and rigorous proofs and methods for many facts on polynomials and formal power series; some of them as sample are treated in this note. Besides the universal properties of the monoid rings and polynomial rings, a universal property for the formal power series rings is also established.
\end{abstract}

\maketitle

\section{Introduction}

Formal power series specially polynomials are ubiquitous in all fields of mathematics and widely applied across the sciences and engineering. Hence, in the abstract setting, it is necessary to have a systematic theory of polynomials available in hand. In \cite[pp. 104-107]{Lang}, the ring of polynomials is introduced very briefly in a systematic way but the details specially the universal property are not investigated. In \cite[Chap 3, \S5]{Hungerford}, this ring is also defined in a satisfactory way but the approach is not sufficiently general. Unfortunately, in the remaining standard algebra text books, this ring is introduced in an unsystematic way. Beside some harmful affects of this approach in the learning, another defect which arises is that then it is not possible to prove many results specially the universal property of the polynomial rings through this non-canonical way. In this note, we plan to fill all these gaps.\\

This material will help to an interesting reader to obtain a better insight into the polynomials and formal power series. It also provides adequate preparation and motivation for reading some advanced books and articles such as \cite{Bruns}, \cite{Gilmer}, \cite{Nishimura} .\\

\section{Monoid rings}

Let $R$ be a ring and let $(M, \ast)$ be a monoid. Let $R[M]$ be the set of all functions $f:M\rightarrow R$ with the finite support, i.e. $\Supp(f)=\{m\in M: f(m)\neq0\}$ is a finite set. We make $R[M]$ into a ring by defining the following operations over it. For $f,g\in R[M]$ the addition is defined as $f+g:M\rightarrow R$ which maps each $m\in M$ into $f(m)+g(m)$. Note that $\Supp(f+g)\subseteq\Supp(f)\cup\Supp(g)$ and hence $f+g\in R[M]$. The multiplication is defined as $f.g:R\rightarrow M$ which maps each $m\in M$ into $$\sum\limits_{{\substack{(m_{1},m_{2})\in M^{2},\\
m_{1}\ast m_{2}=m}}}f(m_{1})g(m_{2}).$$
Note that $E_{m}:=\{(m_{1},m_{2})\in M^{2}: m_{1}\ast m_{2}=m, f(m_{1})g(m_{2})\neq0\}\subseteq\Supp(f)\times\Supp(g)$. Hence $E_{m}$ is a finite set and $(f.g)(m)=\sum\limits_{(m_{1},m_{2})\in E_{m}}f(m_{1})g(m_{2})$. Therefore $f.g$ is actually a function from $M$ into $R$. If $m\in\Supp(f.g)$ then $E_{m}$ is nonempty. Thus the map $m\rightsquigarrow E_{m}$ is an injective map from $\Supp(f.g)$ into the power set of $\Supp(f)\times\Supp(g)$. Therefore $\Supp(f.g)$ is a finite set and hence $f.g\in R[M]$. \\

\begin{lemma}\label{lemma I} The above multiplication is associative. \\
\end{lemma}

{\bf Proof.} Suppose $f,g,h\in R[M]$. For each $m\in M$, consider the equivalence relation $\sim$ on the finite set $E:=\{(m_{1},m_{2},m_{3})\in M^{3}: m_{1}\ast m_{2}\ast m_{3}=m, f(m_{1})g(m_{2})h(m_{3})\neq0\}$ which is defined as $(m_{1},m_{2},m_{3})\sim(n_{1},n_{2},n_{3})$ if $m_{1}\ast m_{2}=n_{1}\ast n_{2}$. Then we recalculate the sum $$c:=\sum\limits_{(m_{1},m_{2},m_{3})\in E}f(m_{1})g(m_{2})h(m_{3})$$ by regarding the equivalent classes which are obtained from this relation. Clearly $$c_{(n_{1},n_{2},n_{3})}:=\sum\limits_{{\substack{(m_{1},m_{2},m_{3})\in E,\\ m_{1}\ast m_{2}=m'}}}f(m_{1})g(m_{2})h(m_{3})$$
is the sum of all terms of $c$ whose indices are belonging to a typical equivalence class with a representation $(n_{1},n_{2},n_{3})$ where $m':=n_{1}\ast n_{2}$. Suppose there exists some $(m_{1},m_{2},m_{3})\in M^{3}\setminus E$ such that $m_{1}\ast m_{2}=m'$, $m'\ast m_{3}=m$ and either $f(m_{1})g(m_{2})\neq0$ or $h(m_{3})\neq0$. Then clearly $f(m_{1})g(m_{2})h(m_{3})=0$. Therefore
$$c_{(n_{1},n_{2},n_{3})}=\sum\limits_{{\substack{m_{3}\in M,\\
m'\ast m_{3}=m}}}\Big(\sum\limits_{{\substack{(m_{1},m_{2})\in M^{2},\\m_{1}\ast m_{2}=m'}}}f(m_{1})g(m_{2})\Big)h(m_{3})=$$$$(f.g)(m')
\Big(\sum\limits_{{\substack{m_{3}\in M,\\ m'\ast m_{3}=m}}}h(m_{3})\Big).$$
Let $\{(n_{1k},n_{2k},n_{3k}): 1\leq k\leq s\}$ be a choice set from the distinct equivalent classes. Let $m'_{k}:=n_{1k}\ast n_{2k}$ and $M':=\{m'_{1},...,m'_{s}\}$. Suppose there is some $(m',m_{3})\in M''\times M$ such that $m'\ast m_{3}=m$ where $M'':=M\setminus M'$. Then clearly $f(m_{1})g(m_{2})h(m_{3})=0$ for each $(m_{1},m_{2})\in M^{2}$ with $m_{1}\ast m_{2}=m'$. Therefore $$\sum\limits_{{\substack{(m',m_{3})\in M''\times M,\\
m'\ast m_{3}=m}}}(f.g)(m')h(m_{3})=0.$$
Now we may write
$$c=\sum\limits_{k=1}^{s}c_{(n_{1k},n_{2k},n_{3k})}=
\sum\limits_{k=1}^{s}(f.g)(m'_{k})\Big(\sum\limits_{{\substack{m_{3}\in M,\\
m_{k}'\ast m_{3}=m}}}h(m_{3})\Big)=\sum\limits_{{\substack{(m',m_{3})\in M'\times M,\\m'\ast m_{3}=m}}}$$$$(f.g)(m')h(m_{3})=
\sum\limits_{{\substack{(m',m_{3})\in M^{2},\\
m'\ast m_{3}=m}}}(f.g)(m')h(m_{3})=\big((f.g).h\big)(m).$$
Similarly, by considering the equivalence relation $\sim'$ over $E$ defined as $(m_{1},m_{2},m_{3})\sim'(n_{1},n_{2},n_{3})$ if $m_{2}\ast m_{3}=n_{2}\ast n_{3}$ then we obtain that $c=\big(f.(g.h)\big)(m)$. $\Box$ \\

The remaining axioms are obvious. Therefore $R[M]$ admits a ring structure by these operations. It is called the monoid ring of $M$ over $R$. In particular, if $G$ is a group then $R[G]$ is called the group ring of $G$ over $R$.
Let $e\in M$ be the identity element of $M$, there is a canonical injective ring map $\eta:R\rightarrow R[M]$ which maps each $a\in R$ into $\eta_{a}:M\rightarrow R$ where $\eta_{a}(e)=a$ and $\eta_{a}(m)=0$ for $m\neq e$. We may identify $R$ as a sub-ring of $R[M]$ by regarding this embedding. If $R$ is a non-zero ring then $\eta$ is an isomorphism if and only if $M$ is trivial. Moreover there is a canonical injective map of monoids $\delta:M\rightarrow R[M]$ to the multiplicative monoid of $R[M]$ which maps each $m\in M$ into $\delta_{m}:M\rightarrow R$ where $\delta_{m}(m')=\delta_{m,m'}$ is the Kronecker delta. The ring $R[M]$ is commutative if and only if $R$ is commutative and $M$ is abelian. Since if
$R[M]$ is commutative then $\delta_{m\ast n}=\delta_{m}.\delta_{n}=\delta_{n}.\delta_{m}=\delta_{n\ast m}$ so $m\ast n=n\ast m$ similarly $\eta_{ab}=\eta_{a}.\eta_{b}=\eta_{b}.\eta_{a}=\eta_{ba}$ and so $ab=ba$.
The converse is also straightforward. Clearly $\eta_{a}.f=af$ for all $a\in R$ and all $f\in R[M]$ where the function $af:M\rightarrow R$ maps each $m\in M$ into $af(m)$.  Moreover $\eta_{a}=a\delta_{e}$ for all $a\in R$. \\

\begin{lemma} The ring $R[M]$ is a free $R-$module with the basis $\{\delta_{m} : m\in M\}$.\\
\end{lemma}

{\bf Proof.} Obviously $R[M]$ has an $R-$module structure through $\eta$.
Each $f\in R[M]$ can be written uniquely as $f=\sum\limits_{m\in M}f(m)\delta_{m}$. $\Box$ \\

\begin{remark}\label{Remark I} Let $\{a_{i}: i\in I\}$ be a set of elements of an additive monoid where $a_{i}$ is the identity element for all but a finite number of indices $i$ and let $\sim$ be an equivalence relation on the index set $I$. Then $$\sum\limits_{i\in I}a_{i}=\sum\limits_{C\in I/\sim}\Big(\sum\limits_{i\in C}a_{i}\Big).$$\\
\end{remark}

Clearly $\eta_{a}.\delta_{m}=\delta_{m}.\eta_{a}$ for all $a\in R$ and all $m\in M$, and the triple $(R[M], \eta, \delta)$ satisfies in the following universal property:\\

\begin{theorem}\label{theorem I} For any such triple $(S, \phi, \psi)$ where $\phi:R\rightarrow S$ is a ring map and $\psi:M\rightarrow S$ is a homomorphism of monoids to the multiplicative monoid of $S$ such that $\phi(a)\psi(m)=\psi(m)\phi(a)$ for all $a\in R$ and all $m\in M$ then there exists a unique ring map $\theta:R[M]\rightarrow S$ such that $\phi=\theta\circ\eta$ and $\psi=\theta\circ\delta$.\\
\end{theorem}

{\bf Proof.} By the above Lemma, there is a unique $R-$linear map $\theta:R[M]\rightarrow S$ such that $\psi=\theta\circ\delta$. We have $\phi=\theta\circ\eta$ since $\theta$ is $R-$linear.
Clearly $\theta$ maps each $f\in R[M]$ into $\sum\limits_{m\in M}\phi\big(f(m)\big)\psi(m)$. Let $f,g\in R[M]$. We have to show that $\theta(f.g)=\theta(f)\theta(g)$. Clearly $$\theta(f.g)=\sum\limits_{m\in M}\Big(\sum\limits_{{\substack{(m_{1},m_{2})\in M^{2},\\m_{1}\ast m_{2}=m}}}\phi\big(f(m_{1})\big)\psi(m_{1})
\phi\big(g(m_{2})\big)\psi(m_{2})\Big)$$ and $$\theta(f)\theta(g)=\Big(\sum\limits_{m_{1}\in M
}\phi\big(f(m_{1})\big)\psi(m_{1})
\Big)\Big(\sum\limits_{m_{2}\in M
}\phi\big(g(m_{2})\big)\psi(m_{2})
\Big)=$$$$\sum\limits_{(m_{1},m_{2})\in M^{2}
}\phi\big(f(m_{1})\big)\psi(m_{1})
\phi\big(g(m_{2})\big)\psi(m_{2}).$$
Consider the equivalence relation $\sim$ on $M^{2}$ which is defined as $(m_{1},m_{2})\sim(n_{1},n_{2})$ if $m_{1}\ast m_{2}=n_{1}\ast n_{2}$.
The map $[(m_{1},m_{2})]\rightsquigarrow m_{1}\ast m_{2}$ is a bijection from the set of distinct equivalent classes onto $M$. Now, using Remark \ref{Remark I}, we get that $$\theta(f)\theta(g)=\sum\limits_{m\in M}\Big(\sum\limits_{{\substack{(m_{1},m_{2})\in M^{2},\\m_{1}\ast m_{2}=m}}}\phi\big(f(m_{1})\big)\psi(m_{1})
\phi\big(g(m_{2})\big)\psi(m_{2})\Big)=\theta(f.g).$$
Hence $\theta$ is a ring homomorphism. Suppose there is another ring map $\theta':R[M]\rightarrow S$ such that $\phi=\theta'\circ\eta$ and $\psi=\theta'\circ\delta$. Then $\theta'$ is $R-$linear since $\phi=\theta'\circ\eta$. It follows that $\theta'=\theta$ since $\psi=\theta'\circ\delta$. $\Box$ \\

For given two monoids $M$ and $N$ consider the monoid $M\times N$ with the component-wise operation. Then one has:\\

\begin{corollary}\label{corollary III} The monoid ring $(R[M])[N]$ is canonically isomorphic to the monoid ring $R[M\times N]$.\\
\end{corollary}

{\bf Proof.} Consider the function $\psi:M\times N\rightarrow (R[M])[N]$ defined by $\psi(m,n)=\eta'(\delta_{m})\delta'_{n}$ where $\eta':R[M]\rightarrow(R[M])[N]$ and $\delta':N\rightarrow (R[M])[N]$ are the canonical maps. The map $\psi$ is a homomorphism of monoids to the multiplicative monoid of $(R[M])[N]$. Then, using Theorem \ref{theorem I}, the desired isomorphism is produced. Some details omitted. $\Box$ \\

Note that in the above Corollary there is no need for the commutativity assumption of neither $R$ and nor monoids.\\

\begin{corollary} For a fixed ring $R$ then the assignment $M\rightsquigarrow R[M]$ is a (covariant) functor from the category of monoids to the category of rings.\\
\end{corollary}

{\bf Proof.} It is an immediate consequence of Theorem \ref{theorem I}. $\Box$ \\

It is worthy to mention that every injective homomorphism of monoids $\psi:M\rightarrow N$ induces a ring map $R[N]\rightarrow R[M]$ given by $f\rightsquigarrow f\circ\psi$. Note that the induced ring map is not necessarily injective.\\

\begin{corollary} For a fixed monoid $M$ then the assignment $R\rightsquigarrow R[M]$ is a (covariant) functor from the category of rings to itself. $\Box$\\
\end{corollary}

\begin{corollary} The map $R[M]\rightarrow R$ given by $f\rightsquigarrow\sum\limits_{m\in M}f(m)$ is a ring homomorphism (it is called the augmentation map).\\
\end{corollary}

{\bf Proof.} Let $\phi:R\rightarrow R$ be the identity map and let $\psi:M\rightarrow R$ be the constant function which maps each element of $M$ into the unit of $R$. Then apply Theorem \ref{theorem I}. $\Box$ \\

\begin{proposition} The kernel of the augmentation map $\theta:R[M]\rightarrow R$ is a free $R-$module with the basis $\{\delta_{m}-1: m\in M, m\neq e\}$.\\
\end{proposition}

{\bf Proof.} Let $A:=R[M]$. If $f\in\Ker\theta$ then $\sum\limits_{m\in M}f(m)=0$. We have then $f=\sum\limits_{m\in M}f(m)\delta_{m}-\big(\sum\limits_{m\in M}f(m)\big)1_{A}=\sum\limits_{m\in M}f(m)\big(\delta_{m}-1_{A}\big)$. Suppose $$\sum\limits_{{\substack{m\in M,\\m\neq e}}}r_{m}(\delta_{m}-1_{A})=0$$ where $r_{m}=0$ for all but a finite number of indices $m$. For each $n\in M$ with $n\neq e$ we have $r_{n}=\big(\sum\limits_{{\substack{m\in M,\\m\neq e}}}r_{m}(\delta_{m}-1_{A})\big)(n)=0$.    $ \Box$\\

\section{Ring of polynomials}

Let $R$ be a ring, $\omega=\{0,1,2,...\}$ the set of natural numbers and $I$ a set. Consider the monoid ring $R[M]$ where
$M$ is the additive monoid consisting of all sequences $(s_{i})_{i\in I}\in\prod\limits_{i\in I}\omega$ such that $s_{i}=0$ for all but a finite number of indices $i$. For each $i\in I$ set $x_{i}:=\delta(\epsilon_{i})$ where $\epsilon_{i}$ is the $i-$th standard unit vector of $M$ and $\delta:M\rightarrow R[M]$ is the canonical map. Each $m=(s_{i})\in M$ can be written as a finite sum $m=\sum\limits_{i\in I}s_{i}\epsilon_{i}$
since $s_{i}=0$ for all but a finite number of indices $i$. We have $\delta(m)=\prod\limits_{i\in I}x^{s_{i}}_{i}$ which we shall denote it simply by $x^{m}$. Therefore each $f\in R[M]$ can be written uniquely as $f=\sum\limits_{m\in M}r_{m}x^{m}$ where $r_{m}:=f(m)$. The ring $R[M]$ is often denoted by $R[x_{i}: i\in I]$ and it is called the ring of polynomials over $R$ with the variables $x_{i}$.\\

Specially, let $I:=\{1,2,3,...,n\}$ for a fixed natural number $n$. Then $M=\omega^{n}$,
$\delta(s_{1},...,s_{n})=x_{1}^{s_{1}}...x^{s_{n}}_{n}$ for all $(s_{1},...,s_{n})\in\omega^{n}$. Also $R[\omega^{n}]$ is denoted by $R[x_{1},...,x_{n}]$. Finally, each $f\in R[x_{1},...,x_{n}]$ can be expressed uniquely as $$f=\sum\limits_{(s_{1},...,s_{n})
\in\omega^{n}}r_{s_{1},...,s_{n}}x_{1}^{s_{1}}...x^{s_{n}}_{n}$$ where $r_{s_{1},...,s_{n}}:=f(s_{1},...,s_{n})$.\\

It is thus seen that in the ring of polynomials, the variables have rigorous mathematical definitions while in many books these are defined informally and as indeterminate objects which are not so pleasant to a mathematician. This insight is just one of the beneficial aspects of the systematic study of a theory that its traces can be found abundantly in Grothendieck's style of mathematics. \\

The following is the universal property of the polynomial rings.\\

\begin{corollary}\label{corollary II} Let $\phi:R\rightarrow A$ be a homomorphism of rings with $A$ commutative and let $\{a_{i}: i\in I\}$ be a set of elements of $A$. Then there exists a unique ring map $\theta:R[x_{i}: i\in I]\rightarrow A$ such that $\phi=\theta\circ\eta$ and $\theta(x_{i})=a_{i}$ for all $i$.\\
\end{corollary}

{\bf Proof.} Let $M$ be as defined in the above, then for each $m=(s_{i})\in M$ there is a finite subset $J$ of $I$ such that $s_{i}=0$ for all $i\in I\setminus J$. Now consider the map $\psi: M\rightarrow A$ which maps each $m$ into $\prod\limits_{i\in J}a_{i}^{s_{i}}$. It is clearly a homomorphism of monoids from the additive monoid of $M$ into the multiplicative monoid of $A$. Hence by Theorem \ref{theorem I}, there is a (unique) ring map $\theta:R[x_{i}: i\in I]\rightarrow A$ such that $\phi=\theta\circ\eta$ and $\psi=\theta\circ\delta$. We have $\theta(x_{i})=\theta\big(\delta(\epsilon_{i})\big)=\psi(\epsilon_{i})=a_{i}$.
Suppose there is another ring map $\theta':R[x_{i}: i\in I]\rightarrow A$
such that $\phi=\theta'\circ\eta$ and $\theta'(x_{i})=a_{i}$ for all $i$. Clearly $\theta'$ is $R-$linear and $\psi=\theta'\circ\delta$. Therefore $\theta'=\theta$.
$\Box$ \\

In Corollary \ref{corollary II}, the image of $\theta$ is a subring of $A$, it is called the $R-$algebra generated by $\{a_{i}: i\in I\}$. A ring map $\phi:R\rightarrow A$ is called of finite type (or, $A$ is a finitely generated $R-$algebra) if there exists a finite set $I$ such that $\theta$ is surjective.\\

If $J\subseteq J'$ then by Corollary \ref{corollary II}, there is a unique ring map $\theta_{J,J'}:R[x_{i}: i\in J]\rightarrow R[y_{j}: j\in J']$ such that $\eta_{J'}=\theta_{J,J'}\circ\eta_{J}$ and $\theta_{J,J'}(x_{i})=y_{i}$ for all $i\in J$ where $\eta_{J}:R\rightarrow R[x_{i}: i\in J]$ is the canonical map.
Let $\Lambda$ be a set of subsets of $I$ such that it is a directed poset by inclusion and $I=\bigcup_{J\in\Lambda}J$ (e.g. $\Lambda$ could be the finite subsets of $I$). The polynomial ring $R[X_{i}: i\in I]$ together with the maps $\lambda_{J}:=\theta_{J,I}$ is the inductive (direct) limit of the inductive system $(R[x_{i}: i\in J], \theta_{J,J'})_{J,J'\in\Lambda}$.
This realization has some applications. For instance, one has:\\

\begin{proposition}\label{prop I} If $R$ is an integral domain then $R[x_{i}: i\in I]$ is as well.\\
\end{proposition}

{\bf Proof.} It is easy to see that the inductive limit of every inductive system of domains is a domain. This reduces the assertion to the finite case. Let $J$ be a finite subset of $I$. In order to prove the assertion for $R[x_{i}: i\in J]$, by induction and Corollary \ref{corollary III}, it suffices to prove it for the one variable polynomial ring $R[x]$. Take two non-zero elements $f,g\in R[x]$. Let $m$ (resp. $n$) be the least natural number such that $f(m)\neq0$ (resp. $g(n)\neq0$). Then $(f.g)(m+n)=f(m)g(n)\neq0$. $\Box$ \\

Let $M$ be as defined in the beginning of this section and let $R[[M]]$ be the set of all functions from $M$ into $R$. This set forms a ring structure by defining the same operations of $R[M]$ over it. It is called the formal power series ring of $R$ with the variables $x_{i}$ and denoted by $R[[x_{i}: i\in I]]$. To see the associativity of the multiplication a similar argument as used in the proof of Lemma \ref{lemma I} is applicable by regarding the fact that for each $m\in M$ then $\{(u,v,w)\in M^{3}: u+v+w=m\}$ is a finite set. Also note that for $f,g\in R[[x_{i}: i\in I]]$ and $u=(u_{i})\in M$ then the summation  $$(f.g)(u)=\sum\limits_
{{\substack{v,w\in M,\\
v+w=u}}}f(v)g(w)$$ is taken over the finite set $$\{(v,w)\in M^{2}: 0\leq v_{i}, w_{i}\leq u_{i}, v_{i}+w_{i}=u_{i}\}$$
where $v=(v_{i})$ and $w=(w_{i})$.
The similar summation is not necessarily definable for a general monoid $M$. Hence it is not always possible to put a ring structure over the set of all functions from a general monoid $M$ into $R$ which admits $R[M]$ as a sub-ring.\\

\begin{remark} For convenience, we may denote each $f\in S:=R[[x_{i}: i\in I]]$ as a sequence $(r_{m})_{m\in M}$ it is also denoted by the notation $\sum\limits_{m\in M}r_{m}x^{m}$ where $r_{m}:=f(m)$ and $x^{m}:=\delta(m)$. The expression $\sum\limits_{m\in M}r_{m}x^{m}$ should not be confused with the usual summation;
it is just a representation of $f$ and in general does not have any mathematical meaning. If $I=\{1,2,..,n\}$ is a finite set then $f$ is also denoted by the notation $f(x_{1},...,x_{n})$. This notation is a little confusing since $f$ is not a function of the variables $x_{1},...,x_{n}$, hence we have to be a little careful when using this notation. Now we observe some examples. In $R[[x]]$ consider the elements $f,g:\omega\rightarrow R$  defined by $f(n):=(-1)^{n}$ and $g(n):=1$ for $n=0,1$ and otherwise $g(n):=0$. By the definition of the multiplication, $f.g=1$. Therefore, in $R[[x]]$, the following identity is obtained:
$$\frac{1}{1+x}=1-x+x^{2}-x^{3}+...,$$
similarly the identity:\\ $$\frac{1}{1-x}=1+x+x^{2}+x^{3}+...$$ holds in $R[[x]]$.
As another type of example, for a fixed $k\in I$ and for a fixed natural number $p\geq0$, we define the function $f_{p,k}:M\rightarrow R$ for each $m=(s_{i})_{i\in I}\in M$ as
$$f_{p,k}(m) :=\frac{(s_{k}+p)!}{(s_{k})!}f(m+p\epsilon_{k}).$$
The map $\Supp(f_{p,k})\rightarrow\Supp(f)$ given by $m\rightsquigarrow m+p\epsilon_{k}$ is injective. Therefore if $f$ is a polynomial then $f_{p,k}$ is as well. The function $f_{p,k}$ is called the formal partial derivative of order $p$ of $f$ w.r.t. the variable $x_{k}$ and  it is denoted by $$\frac{\partial^{p}}{\partial x_{k}}(f)$$ or by $\partial_{k}^{p}f$ or simply by $\partial^{p}f$ if there is no confusion on the variable $x_{k}$. This gives a new general formula to find higher partial derivatives $\partial^{p}f$ directly from original $f$. As specific examples, let $F$ be a field (or a division ring) of characteristic zero, let $a\in F$ and consider the function $\exp_{a}:\omega\rightarrow F$ defined by $\exp_{a}(n):=a^{n}/n!$ then clearly $\partial^{p}\exp_{a}=a^{p}\exp_{a}$. The function $$\exp_{a}=1+ax+\frac{(ax)^{2}}{2!}+\frac{(ax)^{3}}{3!}+...$$ is called the formal exponential function over $F$ (here $s$ means that $s.1_{F}$ for all $s\in\omega$). We may also define the functions $f, g:\omega\rightarrow F$ as $f(n):=(-1)^{s}/n!$ and $g(n):=0$ whenever $n=2s+1$ is odd and $f(n):=0$ and $g(n):=(-1)^{s}/n!$ whenever $n=2s$ is even with $s\in\omega$. Clearly $f(n)=-(n+1)g(n+1)$ and $\partial^{p}f=c_{p}f$ and $\partial^{p}g=c'_{p}g$ whenever $p$ is even and otherwise $\partial^{p}f=c_{p}g$ and $\partial^{p}g=c'_{p}f$ where $c_{p}:=1$ if the remainder of $p$ modulo $4$ is belonging to $\{0,1\}$ and otherwise $c_{p}:=-1$, similarly $c'_{p}:=1$ if the remainder of $p$ modulo $4$ is belonging to $\{0,3\}$ and otherwise $c'_{p}:=-1$. Similarly above $$f=x-\frac{x^{3}}{3!}+\frac{x^{5}}{5!}-\frac{x^{7}}{7!}+...$$ and
$$g=1-\frac{x^{2}}{2!}+\frac{x^{4}}{4!}-\frac{x^{6}}{6!}+...$$
are called the formal sinus and co-sinus functions over $F$, respectively. If there is some element $\lambda\in F$ such that $\lambda^{2}=-1$ (e.g. if $F$ is algebraically closed) then the Euler formula $\lambda f+g=\exp_{\lambda}$ in $F[[x]]$ is easily established. It follows that $$f=\frac{\exp_{\lambda}-\exp_{-\lambda}}{2\lambda},$$
$$g=\frac{\exp_{\lambda}+\exp_{-\lambda}}{2}$$ and so $f^{2}+g^{2}=1$.
It is thus seen that one can develop this theory further and recover formally all of the power series, without no worry about their convergence, which are appearing frequently in mathematical analysis in the form of complex or real valued series. One can also show that the operation $\partial^{p}:S\rightarrow S$, for each $f,g\in S$, has the following properties:\\
$\mathbf{(i)}$ $\partial^{0}$ is the identity map and if $p\geq1$ then $\partial^{p}=\partial^{p-1}\circ\partial$.\\
$\mathbf{(ii)}$ $\partial^{p}$ is additive and $\partial^{p}(f.g)=\sum\limits_{s=0}^{p}\binom{p}{s}
(\partial^{p-s}f).(\partial^{s}g)$. \\
$\mathbf{(iii)}$ $\partial_{k}(x_{\ell})=\delta_{k,\ell}$ and
$\partial_{k}^{p}\circ\partial_{\ell}^{q}=
\partial_{\ell}^{q}\circ\partial_{k}^{p}$ for all $k,\ell\in I$ and all $p,q\in\omega$.\\
\end{remark}

Let $I$ and $J$ be two sets, let $M$ be as defined in the beginning of this section and let $N$ be the set of all sequences $(s_{j})\in\prod\limits_{j\in J}\omega$ such that $s_{j}=0$ for all but a finite number of indices $j$. Then we have:\\

\begin{proposition}\label{corollary IV} The formal power series rings $T:=(R[[M]])[[N]]$ and $S:=R[[M\times N]]$ are isomorphic. \\
\end{proposition}

{\bf Proof.} Consider the function $\psi:T\rightarrow S$ defined as $f\rightsquigarrow f_{\ast}$ where
the function $f_{\ast}:M\times N\rightarrow R$ maps each $\big((s_{i}), (t_{j})\big)\in M\times N$ into $f_{(t_{j})}\big((s_{i})\big)$ and $f_{(t_{j})}:=f\big((t_{j})\big)$. It is not hard to see that $\psi$ is a ring homomorphism and its inverse is built as follows.  Consider the function $\phi:S\rightarrow T$ defined as $f\rightsquigarrow f^{\ast}$ where $f^{\ast}:N\rightarrow R[[M]]$ maps each $(t_{j})\in N$ to the function $f^{\ast}_{(t_{j})}: M\rightarrow R$ which maps each $(s_{i})\in M$ into $f\big((s_{i}), (t_{j})\big)$. $\Box$ \\

\begin{proposition} If $R$ is an integral domain then $R[[x_{1},...,x_{n}]]$ is as well.\\
\end{proposition}

{\bf Proof.} By the induction and Proposition \ref{corollary IV}, it suffices to prove the assertion for $R[[x]]$.  This is proven exactly like Proposition \ref{prop I}. $\Box$ \\

\begin{remark} If $J\subseteq I$ then the (restriction) map $\lambda:M\rightarrow N$ given by $(s_{i})_{i\in I}\rightsquigarrow(s_{i})_{i\in J}$ is a homomorphism of monoids. It is also induces an injective map $R[[N]]\rightarrow R[[M]]$ given by $f\rightsquigarrow f\circ\lambda$. But it is not necessarily a ring homomorphism even if $R$ is commutative.\\
\end{remark}

In the rest of the article $R$ is a commutative ring, $S:=R[[x_{1},...,x_{n}]]$ and $J:=(x_{1},...,x_{n})$ as an ideal of $S$. We have then: \\

\begin{theorem}\label{theorem ivi} The ring $S$ is the completion of $R[x_{1},...,x_{n}]$ w.r.t. the $I-$adic topology where $I:=(x_{1},...,x_{n})$.\\
\end{theorem}

{\bf Proof.} Let $f\in S$. For each natural number $p\geq1$, we define the function
$f_{p}:\omega^{n}\rightarrow R$ as $f_{p}(s_{1},...,s_{n}):=f(s_{1},...,s_{n})$ whenever $\sum\limits_{i=1}^{n}s_{i}< p$ and otherwise $f_{p}(s_{1},...,s_{n}):=0$. Note that $\{(s_{1},...,s_{n})\in\omega^{n}: \sum\limits_{i=1}^{n}s_{i}< p\}$ is a finite set since it is a subset of $\{(s_{1},...,s_{n})\in\omega^{n}: 0\leq s_{i}\leq p-1\}$ which has $p^{n}$ elements. Therefore $f_{p}\in R[x_{1},...,x_{n}]$. In fact $$f_{p}=\sum\limits_{{\substack{(s_{1},...,s_{n})\in\omega^{n},\\
s_{1}+...+s_{n}<p}}}f(s_{1},...,s_{n})x^{s_{1}}_{1}...x^{s_{n}}_{n}.$$
Then consider the map $$\pi_{p}:S\rightarrow R[x_{1},...,x_{n}]/I^{p}$$ defined by $f\rightsquigarrow f_{p}+I^{p}$. It is a ring homomorphism because $(f+g)_{p}=f_{p}+g_{p}$ and $f_{p}g_{p}-(fg)_{p}\in I^{p}$. Moreover  $f_{p+1}-f_{p}\in I^{p}$. Therefore, by the universal property of the projective (inverse) limits, there is a (unique) ring map $$\phi:S\rightarrow\projlim_{p\geq1}R[x_{1},...,x_{n}]/I^{p}$$
which maps each $f$ into $(f_{p}+I^{p})_{p\geq1}$. Clearly it is injective since if $f_{p}\in I^{p}$ then $f_{p}=0$. Take $$(g_{p}+I^{p})_{p\geq1}\in\projlim_{p\geq1}R[x_{1},...,x_{n}]/I^{p}.$$
Then clearly $$g_{p}=\sum\limits_{{\substack{(s_{1},...,s_{n})\in\omega^{n},\\
s_{1}+...+s_{n}<p}}}g_{p}(s_{1},...,s_{n})x^{s_{1}}_{1}...x^{s_{n}}_{n}$$ and $$g_{p+1}=g_{p}+\sum\limits_{{\substack{(s_{1},...,s_{n})\in\omega^{n},\\
s_{1}+...+s_{n}=p}}}g_{p+1}(s_{1},...,s_{n})x^{s_{1}}_{1}...x^{s_{n}}_{n}.$$
Now we define $f:\omega^{n}\rightarrow R$ as $f(s_{1},...,s_{n}):=g_{p+1}(s_{1},...,s_{n})$ where $p:=\sum\limits_{i=1}^{n}s_{i}$.
Note that if $q\geq p+1$ then $g_{q}(s_{1},...,s_{n})=g_{p+1}(s_{1},...,s_{n})$ and so $f_{p}=g_{p}$. Thus $\phi$ maps $f$ into $(g_{p}+I^{p})_{p\geq1}$ and so it is surjective. $\Box$ \\

Let $\mathfrak{a}$ be an ideal of a commutative ring $A$. Recall that the ring $A$ is complete w.r.t. the $\mathfrak{a}-$adic topology if and only if the canonical map $\pi:A\rightarrow\projlim_{p\geq1}A/\mathfrak{a}^{p}$ given by $a\rightsquigarrow(a+\mathfrak{a}^{p})_{p\geq1}$ is bijective. \\

\begin{corollary}\label{th 44} The ring $S$ is complete w.r.t. the $J-$adic topology.\\
\end{corollary}

{\bf Proof.} First note that if $\mathfrak{a}$ is an ideal of a commutative ring $A$ and $B:=\projlim_{p\geq1}A/\mathfrak{a}^{p}$ is the completion of $A$ w.r.t. the $\mathfrak{a}-$adic topology then the ring $B$ is complete w.r.t. the $\mathfrak{a}B-$adic topology where $\mathfrak{a}B$ is the extension of $\mathfrak{a}$ under the canonical map $\pi:A\rightarrow B$. Indeed, this fact holds more generally, e.g. see \cite[Proposition 10.3]{Atiyah}. Now the assertion is clear from this fact and Theorem \ref{theorem ivi}. $\Box$ \\

Finding a direct proof to the following subtle result would be certainly a serious challenge for the readers.\\

\begin{corollary}\label{coro 45} If $f\in S$ then $f-f_{p}\in J^{p}$ for all $p\geq1$.\\
\end{corollary}

{\bf Proof.} By the above Corollary, there is some $g\in S$ such that $g-f_{p}\in J^{p}$ for all $p\geq1$. We show that $g=f$. Take $(s_{1},...,s_{n})\in\omega^{n}$ and let $q:=s_{1}+...+s_{n}$. We have $g=f_{q+1}+\sum\limits_{\substack{(t_{1},...,t_{n})\in\omega^{n},\\
t_{1}+...+t_{n}=q+1}}h_{t_{1},...,t_{n}}.(x^{t_{1}}_{1}...x^{t_{n}}_{n})$ where $h_{t_{1},...,t_{n}}\in S$. It follows that $g(s_{1},...,s_{n})=f(s_{1},...,s_{n})$. $\Box$ \\

\begin{corollary} The map $R[x_{1},...,x_{n}]/I^{p}\rightarrow S/J^{p}$ induced by the canonical injection $i:R[x_{1},...,x_{n}]\rightarrow S$ is an isomorphism.\\
\end{corollary}

{\bf Proof.} It is obvious from the above Corollary. $\Box$ \\

The triple $(S,\eta',J)$ satisfies in the following universal property where $\eta':=i\circ\eta$ and $\eta:R\rightarrow R[x_{1},...,x_{n}]$ is the canonical map.\\

\begin{theorem}\label{Th 20} For each such triple $(A,\phi,\mathfrak{a})$, where $\phi:R\rightarrow A$ is a ring homomorphism into the commutative ring $A$ and $\mathfrak{a}$ is an ideal of $A$ such that the $\mathfrak{a}-$adic topology over it is complete, and for each finite subset $\{a_{1},...,a_{n}\}$ of $\mathfrak{a}$ then there exists a unique ring map $\psi:S\rightarrow A$ such that $\phi=\psi\circ\eta'$ and $\psi(x_{i})=a_{i}$ for all $i$. \\
\end{theorem}

{\bf Proof.} By Corollary \ref{corollary II}, there is a (unique) ring map $$\theta:R[x_{1},...,x_{n}]\rightarrow A$$ such that $\phi=\theta\circ\eta$ and $\theta(x_{i})=a_{i}$ for all $i$. Then, by the universal property of the projective limits, there is a unique ring map $\psi':S\rightarrow T:=\projlim_{p\geq1}A/\mathfrak{a}^{p}$ such that $\theta_{p}\circ\pi_{p}=\pi'_{p}\circ\psi'$ for all $p\geq1$ where $\pi'_{p}:T\rightarrow A/\mathfrak{a}^{p}$ and $\pi_{p}:S\rightarrow R[x_{1},...,x_{n}]/I^{p}$ are the canonical projections for the latter see Theorem \ref{theorem ivi} and $\theta_{p}:R[x_{1},...,x_{n}]/I^{p}\rightarrow A/\mathfrak{a}^{p}$ is induced by $\theta$. Clearly $$\pi'_{p}\big(\psi'(f)\big)=
\sum\limits_{\substack{(s_{1},...,s_{n})\in\omega^{n},\\
s_{1}+...+s_{n}<p}}
\phi(r_{s_{1},...,s_{n}})a^{s_{1}}_{1}...a^{s_{n}}_{n}+\mathfrak{a}^{p}$$ for all $f\in S$ and all $p\geq1$ where $r_{s_{1},...,s_{n}}:=f(s_{1},...,s_{n})$. The canonical ring map $\pi:A\rightarrow T$ given by $a\rightsquigarrow (a+\mathfrak{a}^{p})_{p\geq1}$ is bijective since $A$ is $\mathfrak{a}-$complete. Let $\psi:=\pi^{-1}\circ\psi'$. Then it is easy to see that $\psi(x_{i})=a_{i}$ for all $i$ and $\psi\circ\eta'=\phi$. Now suppose there is another ring map $h:S\rightarrow A$ such that $\phi=h\circ\eta'$ and $h(x_{i})=a_{i}$ for all $i$. We have $h(f_{p})=\theta(f_{p})$ and by Corollary \ref{coro 45}, $h(f-f_{p})\in \mathfrak{a}^{p}$ for all $f\in S$. This shows that $\theta_{p}\circ\pi_{p}=\pi'_{p}\circ(\pi\circ h)$ and so $h=\psi$. $\Box$ \\

The above Theorem is called the universal property of the formal power series rings. One can also show that $\psi$ is continuous when $S$ and $A$ are considered as topological spaces via the $J-$adic and $\mathfrak{a}-$adic topologies, respectively.\\

Note that in general, for a given $a\in R$, there is no ring map from $R[[x]]$ to
$R$ which maps $x$ into $a$. In fact there are many facts on $R[x_{1},...,x_{n}]$ which can not be extended to $R[[x_{1},...,x_{n}]]$; roughly speaking many of these failures are due to the danger of the divergence in infinite sums.\\

\begin{remark} One can directly show that the canonical map $$\pi:T\rightarrow\projlim_{p\geq1}T/L^{p}$$ is injective where $T:=R[[x_{i}: i\in I]]$ and
$L:=(x_{i}: i\in I)$ (note that here $I$ is not necessarily a finite set).
Because assume that $f\in\Ker\pi=\bigcap\limits_{p\geq1}L^{p}$. Take $\textbf{s}=(s_{i})_{i\in I}\in M$ and let $q:=\sum\limits_{i\in I}s_{i}$ where $M$ is as defined in the beginning of this section. We have $f\in L^{q+1}$. Hence there are a finitely many elements $g_{k}\in T$ and $h_{k}\in J^{q+1}$ such that $f=\sum\limits_{k}g_{k}.h_{k}$. Note that corresponding with each $h_{k}$ there exists some $(t_{i})_{i\in I}\in M$ such that $\sum\limits_{i\in I}t_{i}=q+1$ and $h_{k}=\prod\limits_{i\in I}x^{t_{i}}_{i}$.
It follows that $(g_{k}.h_{k})(\textbf{s})=0$ and so
$f(\textbf{s})=0$. It is natural to ask, is $\pi$ a surjective map? This is not known for the author. Note that if $f\in T$ and $I$ is an infinite set then $f_{p}$ is not necessarily a polynomial since the set $E_{p}:=\{m=(s_{i})\in M : \sum\limits_{i\in I}s_{i}<p\}$ for $p\geq2$ is an infinite set, e.g. $(p-1)\epsilon_{i}\in E_{p}$ for all $i\in I$. Thus the finiteness of the number of variables in Theorem \ref{theorem ivi} is a crucial point and hence we are not able to prove this result and also Theorem \ref{Th 20} for an arbitrary number of variables. Despite of these, proposing the following question seems reasonable.\\
\end{remark}

$\mathbf{Question.}$ Suppose $\alpha$ is an ordinal number and $R[[x_{i}: i\in\beta]]$ is complete w.r.t. the $(x_{i}: i\in\beta)-$adic topology for all ordinals $\beta$ with $\beta<\alpha$. Is $R[[x_{i}: i\in\alpha]]$ complete w.r.t. the $(x_{i}: i\in\alpha)-$adic topology?\\

Let $I$ be a set. If the above question is true then the transfinite induction would imply that $R[[x_{i}: i\in I]]$ is complete w.r.t. the $(x_{i}: i\in I)-$adic topology.\\

\end{document}